\documentclass{svmult-ddm}

\usepackage{mathptmx}       
\usepackage{helvet}         
\usepackage{courier}        
\usepackage{type1cm}        
\usepackage{graphicx}        

\usepackage[bottom]{footmisc}

\usepackage{amssymb}
\usepackage{amsfonts}
\usepackage{amsbsy}
\usepackage{amscd}
\usepackage{amstext}
\usepackage{dsfont}
\usepackage[english]{babel}
\usepackage{color}
\usepackage{graphics}
\usepackage{epsfig}
\usepackage{subfigure}
\usepackage{wrapfig}
\usepackage{psfrag}
\usepackage{color}
\usepackage{url}
\usepackage{verbatim}
\usepackage{algorithm}
\usepackage{algorithmic}

\usepackage{amsmath}

\begin{document}

\title*{Generating Equidistributed Meshes in 2D via Domain Decomposition}
\titlerunning{2D Meshes via DD}
\author{Ronald D. Haynes\inst{1} \and
Alexander J. M. Howse\inst{1}}
\authorrunning{R. Haynes \and A. Howse}
\institute{\inst{1}Memorial University of Newfoundland, Department of Mathematics \& Statistics, St. John's, NL,  Canada  A1C 5S7
\email{{rhaynes}{z37ajmh}@mun.ca} }
%
%
\maketitle

\abstract*{In this paper we consider Schwarz domain decomposition applied to the generation of 2D spatial meshes by a local 
  equidistribution principle. We briefly review the derivation of the local equidistribution principle and the appropriate choice of 
  boundary conditions. We then introduce classical and optimized Schwarz domain decomposition methods to solve the resulting system of 
  nonlinear equations. The implementation of these iterations are discussed, and we conclude with numerical examples to illustrate the 
  performance of the approach.
}

\section{Introduction}

There are many occasions when the use of a uniform spatial grid would be prohibitively expensive for the numerical solution of partial differential equations (PDEs). In such situations, a popular strategy is to generate an adaptive 
mesh by either varying the number of mesh points, the order of the numerical method, or the location of mesh points throughout the domain,
in order to best resolve the solution. It is the latter of these options, known as \emph{moving mesh methods}, which is our focus. In this case 
the physical PDE of interest is coupled with equations which adjust the position of mesh points to best ``equidistribute'' a particular measure of numerical error. This coupled system of equations is solved to generate the solution and the corresponding mesh simultaneously, see \cite{haynesr_mini_17_huangbook} for a recent overview.

A simple method for adaptive grid generation in two spatial dimensions is outlined in \cite{haynesr_mini_17_2dmesh} by Huang and Sloan, in which a finite difference two dimensional adaptive mesh method is developed by applying a variation of de Boor's equidistribution principle (EP) \cite{haynesr_mini_17_deBoor1, haynesr_mini_17_deBoor2}. The equidistribution principle states that an appropriately chosen mesh should equally distribute some measure of the solution variation or computational error over the entire domain. Mackenzie \cite{haynesr_mini_17_mack} extends upon the work of \cite{haynesr_mini_17_2dmesh} by presenting a finite volume discretization of the mesh equations, as well as an efficient iterative approach for solving these equations, referred to as ``an alternating line Gauss-Seidel relaxation approach''.

In this paper, we propose a parallel domain decomposition (DD) solution of the 2D adaptive method of \cite{haynesr_mini_17_2dmesh}. In Section \ref{haynesr_mini_17_sec:theory} we review the derivation of the mesh PDEs of \cite{haynesr_mini_17_2dmesh} and discuss possible 
boundary conditions. In Sections \ref{haynesr_mini_17_sec:ddtheory} and \ref{haynesr_mini_17_optimized} we present classical and optimized Schwarz methods for the generation of 2D equidistributed meshes, and in Section \ref{haynesr_mini_17_sec:dd:implementation} we describe the numerical implementation of this approach and provide numerical results.

\section{2D Mesh Generation}
\label{haynesr_mini_17_sec:theory}

To begin, we review the derivation of the equations which govern mesh equidistribution in two 
spatial dimensions from \cite{haynesr_mini_17_2dmesh}, defining a mesh in the physical variables $(x,y)$ which \emph{best} resolves a given function $u(x,y)$. Let $\mathbf{x} = [x,y]^T$ be the spatial coordinates of a mesh in a 2D physical domain, $\Omega_p$. We introduce the coordinate transformation $\mathbf{x} = \mathbf{x}(\pmb{\xi})$, where $\pmb{\xi} = [\xi, \eta]^T$ denotes the spatial coordinates on the computational domain, $\Omega_c = [0,1]\times[0,1]$. Here we determine a mesh which equidistributes the arc-length of $u(x,y)$ over $\Omega_p$. The scaled arc-length variation of $u$ along the arc element from $\mathbf{x}$ to $\mathbf{x} + d\mathbf{x}$ can be expressed as
\begin{equation}
ds = [a^2(du)^2 + d\mathbf{x}^T d\mathbf{x}]^{1/2} = [d\mathbf{x}^T \pmb{M} d\mathbf{x}]^{1/2},
\label{haynesr_mini_17_arclength2d}
\end{equation}
where $\pmb{M}  = a^2\pmb{\nabla} u\,\cdot \pmb{\nabla} u^T + \pmb{I}$ and $a\in[0,1]$ is a user specified relaxation parameter. The extreme cases are $a=0$, which produces a uniform mesh, and $a=1$, which produces a mesh equidistributing the arc-length monitor function. Making use of the mesh transformation $\mathbf{x} = \mathbf{x}(\pmb{\xi})$, (\ref{haynesr_mini_17_arclength2d}) can be expressed as
\begin{equation}
ds = [d\pmb{\xi}^T\,\pmb{J}^T\,\pmb{M}\,\pmb{J}\,d\pmb{\xi}]^{1/2},
\label{haynesr_mini_17_arclength2dtransformed}
\end{equation}
where $\pmb{J}$ is the Jacobian of the transformation.

The equidistribution principle follows from (\ref{haynesr_mini_17_arclength2dtransformed}): if $u(\mathbf{x}(\pmb{\xi}))$ is to have the same value $ds$ along any arc element in the computational domain with fixed length $[d\pmb{\xi}^T d\pmb{\xi}]^{1/2}$, then (\ref{haynesr_mini_17_arclength2dtransformed}) must be independent of the coordinate $\pmb{\xi}$. This implies that $\pmb{J}^T\pmb{M}\pmb{J}$ should be independent of $\pmb{\xi}$, or
\begin{equation}
[d\pmb{\xi}^T\,\pmb{J}^T\,\pmb{M}\,\pmb{J}\,d\pmb{\xi}]^{1/2} = [d\pmb{\xi}^T\,\pmb{\widetilde{M}}\,d\pmb{\xi}]^{1/2},
\label{haynesr_mini_17_2dequidist}
\end{equation}
where $\pmb{\widetilde{M}}$ is a constant and hence $\pmb{\xi}$-independent matrix.
If a coordinate transformation can be found which satisfies (\ref{haynesr_mini_17_2dequidist}), $u$ will have the same variation at any point in $\Omega_p$ along any arc of length
$$
\left[\left(\textstyle\frac{\partial\mathbf{x}}{\partial \xi} d\xi + \textstyle\frac{\partial\mathbf{x}}{\partial \eta} d\eta \right)^T \left(\textstyle\frac{\partial\mathbf{x}}{\partial \xi} d\xi + \textstyle\frac{\partial\mathbf{x}}{\partial \eta} d\eta \right)\right]^{1/2}.
$$
A transformation satisfying (\ref{haynesr_mini_17_2dequidist}) for some matrix $\pmb{\widetilde{M}}$ will be called an \emph{equidistribution}, and (\ref{haynesr_mini_17_2dequidist}) an \emph{equidistribution principle}.

Usually (\ref{haynesr_mini_17_2dequidist}) cannot be satisfied by the coordinate transformation on the whole computational domain. However, if (\ref{haynesr_mini_17_2dequidist}) is weakened so that the transformation is only required to satisfy (\ref{haynesr_mini_17_2dequidist}) locally; that is, we only require $\pmb{\widetilde{M}}$ to be constant along a given coordinate line, it is possible to find a local equidistribution on $\Omega_p$. In 2D this leads to the system:
\begin{equation}
\left[\left(
 \begin{array}{c}
  \textstyle\frac{\partial x}{\partial \xi} \\
  \textstyle\frac{\partial y}{\partial \xi}  
 \end{array}
 \right)^T
 \pmb{M}
 \left(
 \begin{array}{c}
  \textstyle\frac{\partial x}{\partial \xi} \\
  \textstyle\frac{\partial y}{\partial \xi}  
 \end{array}
 \right)
 \right]^{1/2} = c_1(\eta), 
\quad\quad
 \left[\left(
 \begin{array}{c}
  \textstyle\frac{\partial x}{\partial \eta} \\
  \textstyle\frac{\partial y}{\partial \eta}  
 \end{array}
 \right)^T
\pmb{M}
 \left(
 \begin{array}{c}
  \textstyle\frac{\partial x}{\partial \eta} \\
  \textstyle\frac{\partial y}{\partial \eta}  
 \end{array}
 \right)
 \right]^{1/2} = c_2(\xi),
 \label{haynesr_mini_17_2dsystem}
\end{equation}
where $c_1(\eta)$ and $c_2(\xi)$ are constant in the $\xi$ and $\eta$ directions respectively.
These constants are eliminated by numerical differencing.

Instead of using the scaled arc-length matrix $\pmb{M}$, in practice we modify $\pmb{M}$ as
$\pmb{M} = k\pmb{\nabla} u \cdot \pmb{\nabla} u^T+ \pmb{I}$, where $k=a^2/(1+b\pmb{\nabla} u^T\pmb{\nabla} u)$.
The parameter $b \geq 0$ is used to prevent problems where extremely small mesh spacing or mesh tangling could occur, that is when $|\pmb{\nabla} u|$ is very large. 

System (\ref{haynesr_mini_17_2dsystem}) will determine the internal mesh points. 
In \cite{haynesr_mini_17_2dmesh} a combination of Dirichlet and Neumann conditions are used along $\partial\Omega_c$:
\begin{eqnarray}
\label{haynesr_mini_17_boundaryconds1}
x(0,\eta) = y(\xi,0) = 0, \quad x(1,\eta) = y(\xi,1) = 1, \\
\label{haynesr_mini_17_boundaryconds2}
\displaystyle\frac{\partial x}{\partial \eta}(\xi,0) = \displaystyle\frac{\partial x}{\partial \eta}(\xi,1) = \displaystyle\frac{\partial y}{\partial \xi}(0,\eta) = \displaystyle\frac{\partial y}{\partial \xi}(1,\eta) = 0,
\end{eqnarray}
where $\xi,\eta \in [0,1]$. The Dirichlet conditions are consistent with the requirement that there are mesh points on the boundary of the domain. The Neumann orthogonality conditions are arbitrary, and in fact can cause smoothness issues near the domain boundaries. As an alternative, we 
follow \cite{haynesr_mini_17_mack} and apply the 1D EP,
\begin{equation}
\label{haynesr_mini_17_1dequidist}
(M(x)x_\xi)_\xi = 0, \quad x(0) = 0, \quad x(1) = 1,
\end{equation}
to determine $x(\xi,0)$, $x(\xi,1)$, $y(0,\eta)$ and $y(1,\eta)$. The 1D analog of the system (\ref{haynesr_mini_17_2dsystem}), given in (\ref{haynesr_mini_17_1dequidist}), has previously been solved by DD methods in  \cite{haynesr_mini_17_SINUM, haynesr_mini_17_haynes12, haynesr_mini_17_ETNA}.

\section{Classical Schwarz Domain Decomposition}
\label{haynesr_mini_17_sec:ddtheory}

\begin{wrapfigure}{R}{0.5\textwidth}
\vspace{-30pt}
  \begin{center}
    \includegraphics[width=0.48\textwidth]{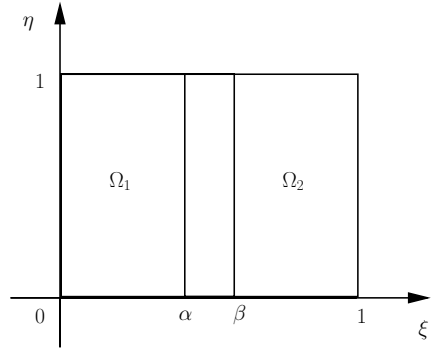}
  \end{center}
  \vspace{-10pt}
  \caption{DD in $\xi$ using in 2 subdomains.}
  \label{haynesr_mini_17_fig1}
  \vspace{-20pt}
\end{wrapfigure}

For the two dimensional mesh adaptation problem, the computational domain $\Omega_c = [0,1]\times[0,1]$, can either be decomposed in just the
$\xi$ or just the $\eta$ directions, or in both directions.  This results in  ``strip'' or ``block'' configurations of 
subdomains respectively. Here we discuss DD applied in the $\xi$ direction only. That is, we decompose the $\xi$ interval $[0,1]$ into subintervals $[\alpha_{\xi}^{i}, \beta_{\xi}^{i}]$, $i=1,\ldots,S$, where $\alpha_{\xi}^{1} = 0$, $\beta_{\xi}^{S} = 1$,
and we assume the subintervals satisfy the overlap conditions:
$$
\alpha_{\xi}^{i} < \alpha_{\xi}^{i+1} < \beta_{\xi}^{i} < \beta_{\xi}^{i+1}.
$$
The resulting decomposition has $S$ subdomains, denoted by $\Omega_i = [\alpha_{\xi}^{i}, \beta_{\xi}^{i}]\times[0,1]$ for $i=1,\ldots,S$. 
The boundary 
conditions (\ref{haynesr_mini_17_boundaryconds1}--\ref{haynesr_mini_17_boundaryconds2}) or (\ref{haynesr_mini_17_1dequidist}) are used 
along the ends of each strip and transmission conditions are specified along the newly created interfaces.

Consider the 2D adaptive mesh system, (\ref{haynesr_mini_17_2dsystem}), for the $S=2$ case.  We split $\Omega_c$  
into subdomains $\Omega_1$ and $\Omega_2$ as in Figure \ref{haynesr_mini_17_fig1}.  Let $x_i^n$ denote the subdomain
solution on $\Omega_i$, for $i=1,2$.  We consider the following DD iteration: for $n=1,2,\ldots$, 
solve
\begin{eqnarray}
\label{haynesr_mini_17_iter_eq1}
\left[\left(
 \begin{array}{c}
  \textstyle\frac{\partial x_i^n}{\partial \xi} \vspace{2pt}\\
  \textstyle\frac{\partial y_i^n}{\partial \xi}  
 \end{array}
 \right)^T
 \pmb{M}(x_i^n, y_i^n)
 \left(
 \begin{array}{c}
  \textstyle\frac{\partial x_i^n}{\partial \xi} \vspace{2pt}\\
  \textstyle\frac{\partial y_i^n}{\partial \xi}  
 \end{array}
 \right)
 \right]^{1/2} = c_1(\eta), 
\\
\label{haynesr_mini_17_iter_eq2}
 \left[\left(
 \begin{array}{c}
  \textstyle\frac{\partial x_i^n}{\partial \eta} \vspace{2pt}\\
  \textstyle\frac{\partial y_i^n}{\partial \eta}  
 \end{array}
 \right)^T
\pmb{M}(x_i^n, y_i^n)
 \left(
 \begin{array}{c}
  \textstyle\frac{\partial x_i^n}{\partial \eta} \vspace{2pt}\\
  \textstyle\frac{\partial y_i^n}{\partial \eta}  
 \end{array}
 \right)
 \right]^{1/2} = c_2(\xi),
\end{eqnarray}
for $i=1,2$ and $\pmb{\xi} \in \Omega_i$. The classical Schwarz method uses the transmission conditions
\begin{eqnarray}
x_1^n(\beta,\eta) &=& x_2^{n-1}(\beta,\eta),\quad y_1^n(\beta,\eta) = y_2^{n-1}(\beta,\eta),\label{haynesr_mini_17_2dsystem:sd1} \\
x_2^n(\alpha,\eta) &=& x_1^{n-1}(\alpha,\eta),\quad y_2^n(\alpha,\eta) = y_1^{n-1}(\alpha,\eta). \label{haynesr_mini_17_2dsystem:sd2}
\end{eqnarray}
On $\partial(\Omega_c \cap \Omega_i)$ the boundary conditions (\ref{haynesr_mini_17_boundaryconds1}) are used, along with the 1D EP 
to determine $x(\xi,0)$, $x(\xi,1)$, $y(0,\eta)$ and $y(1,\eta)$.


Each DD iteration requires a pair of PDEs to be solved, each a ``smaller'' version of the local EP (\ref{haynesr_mini_17_2dsystem}). 
These problems are solved in an iterative manner: given initial approximations to be used along interfaces, 
the PDEs (\ref{haynesr_mini_17_iter_eq1}--\ref{haynesr_mini_17_iter_eq2}) are solved, and then solution information along the interfaces 
is exchanged between subdomains. The PDEs are then solved again, 
now with updated boundary data, and the process repeats. By iterating, 
the subdomain solutions converge to the 
desired solution $\mathbf{x}$ on their respective subdomains. 
As is well known, classical Schwarz requires the subdomains to overlap \cite{haynesr_mini_17_ddhistory}. 

\section{Optimized Boundary Conditions}
\label{haynesr_mini_17_optimized}

Classical Schwarz is known to converge slowly. As a way to remedy this, we propose the use of higher order, Robin type, transmission conditions along the artificial interfaces.  As before, we decompose $\Omega_c = [0,1]\times[0,1]$ into subdomains $\Omega_1 = [0,\beta]\times[0,1]$ and $\Omega_2 = [\alpha, 1]\times[0,1]$, where $\alpha \leq \beta$. 

We define, for any differentiable functions $x(\xi,\eta)$ and $y(\xi,\eta)$, the operators
\begin{eqnarray*}
B_1(x) = x_\xi + p x, & \quad B_2(x) = x_\xi - p x, \\
B_3(x, y) = S_1(x, y) + p x, & \quad B_4(x, y) = S_1(x, y) - p x,
\end{eqnarray*}
where
$$
S_1(x,y) = \sqrt{
\left(\begin{array}{c}
x_{\xi} \\
y_{\xi}
\end{array}\right)^T
M
\left(\begin{array}{c}
x_{\xi} \\
y_{\xi}
\end{array}\right)
},
\quad
M = \frac{a^2 w \cdot w^T}{1 + b w^T \cdot w} + I
$$
and
$$
w = \frac{1}{x_\xi y_\eta - x_\eta y_\xi}\left[u_\xi y_\eta - u_\eta y_\xi,  -u_\xi x_\eta + u_\eta x_\xi \right]^T.
$$
We propose two possible sets of transmission conditions. The first are simple linear Robin conditions using the derivative normal to the artificial boundaries:
\begin{equation}
\label{haynesr_mini_17_4transconds}
\begin{split}
& B_1(x_{1}^n(\beta, \eta)) = B_1(x_{2}^{n-1}(\beta, \eta)), \quad B_1(y_{1}^n(\beta, \eta)) = B_1(y_{2}^{n-1}(\beta, \eta)) \\
& B_2(x_{2}^n(\alpha, \eta)) = B_2(x_{1}^{n-1}(\alpha, \eta)), \quad B_2(y_{2}^n(\alpha, \eta)) = B_2(y_{1}^{n-1}(\alpha, \eta)).
\end{split}
\end{equation}

The second set are of nonlinear Robin type, similar to those used in 
an optimized Schwarz algorithm for 1D mesh generation in 
\cite{haynesr_mini_17_SINUM}.
We replace the $x$ equations of (\ref{haynesr_mini_17_4transconds}) by
\begin{equation}
\label{haynesr_mini_17_iter2:1}
\begin{split}
& B_3(x_{1}^n(\beta, \eta),y_{1}^n(\beta, \eta)) = B_3(x_{2}^{n-1}(\beta, \eta)y_{2}^{n-1}(\beta, \eta)) \\
& B_4(x_{2}^n(\alpha, \eta),y_{2}^n(\alpha, \eta)) = B_4(x_{1}^{n-1}(\alpha, \eta),y_{1}^{n-1}(\alpha, \eta)).
\end{split}
\end{equation}
Note, the mesh PDE (\ref{haynesr_mini_17_iter_eq1}) indicates that the nonlinear term $S_1$ in the operator $B_3$ is constant across the $\xi = \alpha$ and $\xi = \beta$ interfaces. Furthermore, as the system of equations resulting from (\ref{haynesr_mini_17_iter_eq1}-\ref{haynesr_mini_17_iter_eq2}) are already nonlinear, the nonlinear transmission conditions will not have a large impact on the cost of solving the system.

\section{Numerical Implementation and Results}
\label{haynesr_mini_17_sec:dd:implementation}

The local EP (\ref{haynesr_mini_17_2dsystem}), the physical boundary conditions on $\Omega_c$, and the transmission conditions (\ref{haynesr_mini_17_2dsystem:sd1}, \ref{haynesr_mini_17_2dsystem:sd2}), (\ref{haynesr_mini_17_4transconds}) or (\ref{haynesr_mini_17_iter2:1}), are discretized using standard finite differences on a uniform grid in the computational $(\xi, \eta)$ variables. Second order centered differences are used, using the ghost value technique as needed at the boundaries to ensure the scheme is second order. 
The nonlinear transmission conditions require nonlinear, rather than linear, equations to be solved at the interface.  This is not onerous
as the whole system is solved with a Newton iteration.

%

\begin{figure}[h]
\begin{minipage}[b]{0.47\linewidth}
\centering
 \includegraphics[width=.95\linewidth]{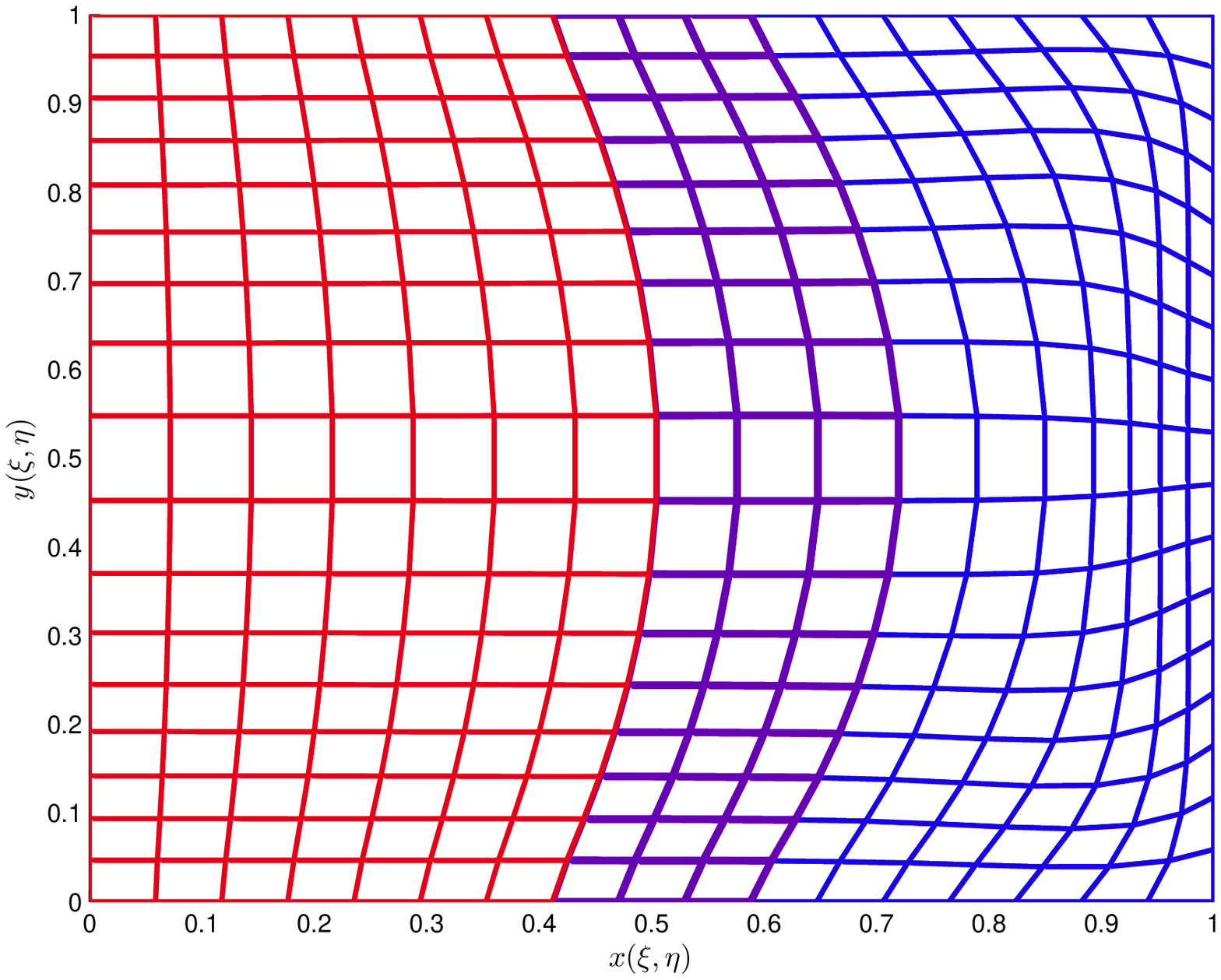}
 \caption{Adaptive mesh generated for the test function.}
\label{haynesr_mini_17_figure_u_view1}
\end{minipage}
\hspace{0.5cm}
\begin{minipage}[b]{0.47\linewidth}
\centering
 \includegraphics[width=.95\linewidth]{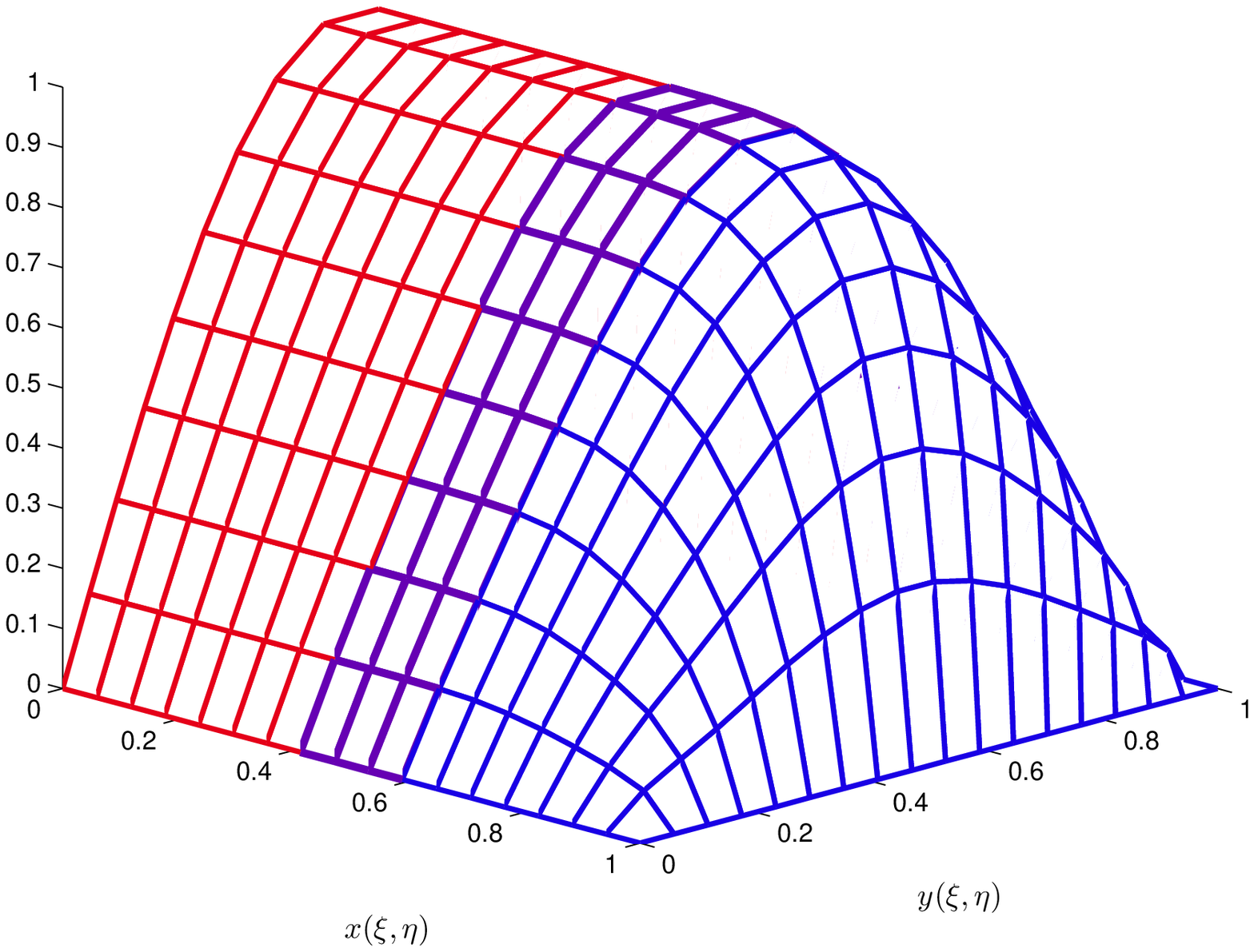}
\caption{The test function plotted using an adaptive mesh.}
\label{haynesr_mini_17_figure_u_view2}
\end{minipage}
\end{figure}

In the examples we use the test function
$u(x,y) =\left[1 - e^{15(x-1)} \right]\sin(\pi y).$
The function is shown, along with its locally equidistributed mesh, in Figures \ref{haynesr_mini_17_figure_u_view1} and \ref{haynesr_mini_17_figure_u_view2}.  
The physical mesh $(x,y)$ is generated by solving (\ref{haynesr_mini_17_2dsystem}) using a grid of $18\times 18$ uniformly spaced 
mesh points in $\Omega_c$. 
For this example, we use an optimized Schwarz iteration, with
transmission conditions (\ref{haynesr_mini_17_4transconds}), on 2 subdomains with 4 points of overlap 
in the $\xi$ direction.  Here the number of points of overlap refers to the number of shared grid points, the overlap
width is approximately half of this number times $\Delta\xi$.
We choose the parameters $a=0.7$, $b = 0.05$ and $p=2.3$. 
The mesh on subdomain 1 is shown in red, on subdomain 2 in blue, and the overlap region in purple. 
In general, the meshes obtained by the different methods will be visually indistinguishable from one another at convergence. To compare the DD methods we will plot their convergence histories.

\begin{figure}[h]
\begin{minipage}[b]{0.47\linewidth}
\centering
\includegraphics[width=.95\linewidth]{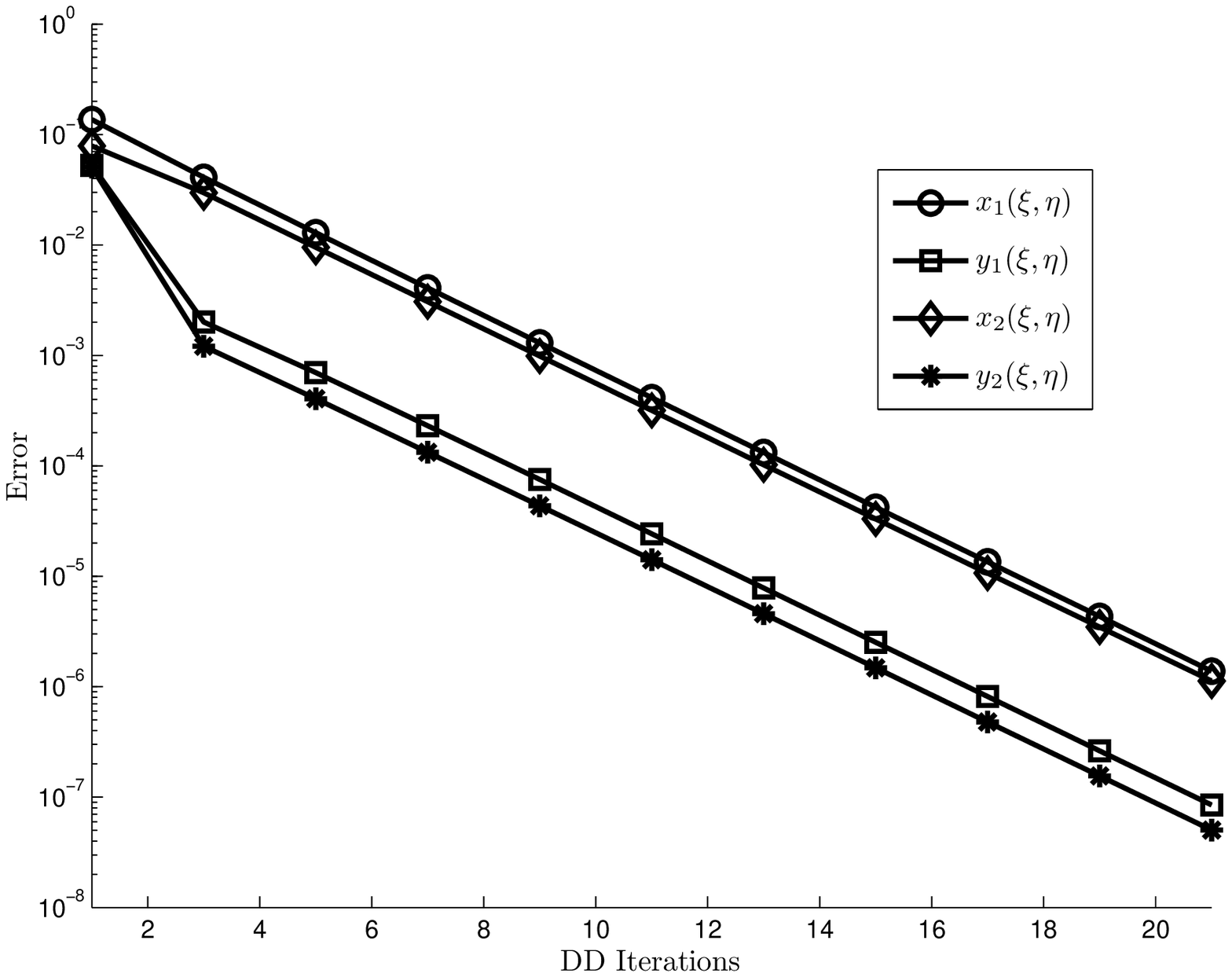}
\caption{Classical Schwarz convergence histories for each part of the solution, $x_{1,2}^n$ and $y_{1,2}^n$.}
\label{haynesr_mini_17_figure_DD_2D_comp1_1}
\end{minipage}
\hspace{0.5cm}
\begin{minipage}[b]{0.47\linewidth}
\centering
 \includegraphics[width=.95\linewidth]{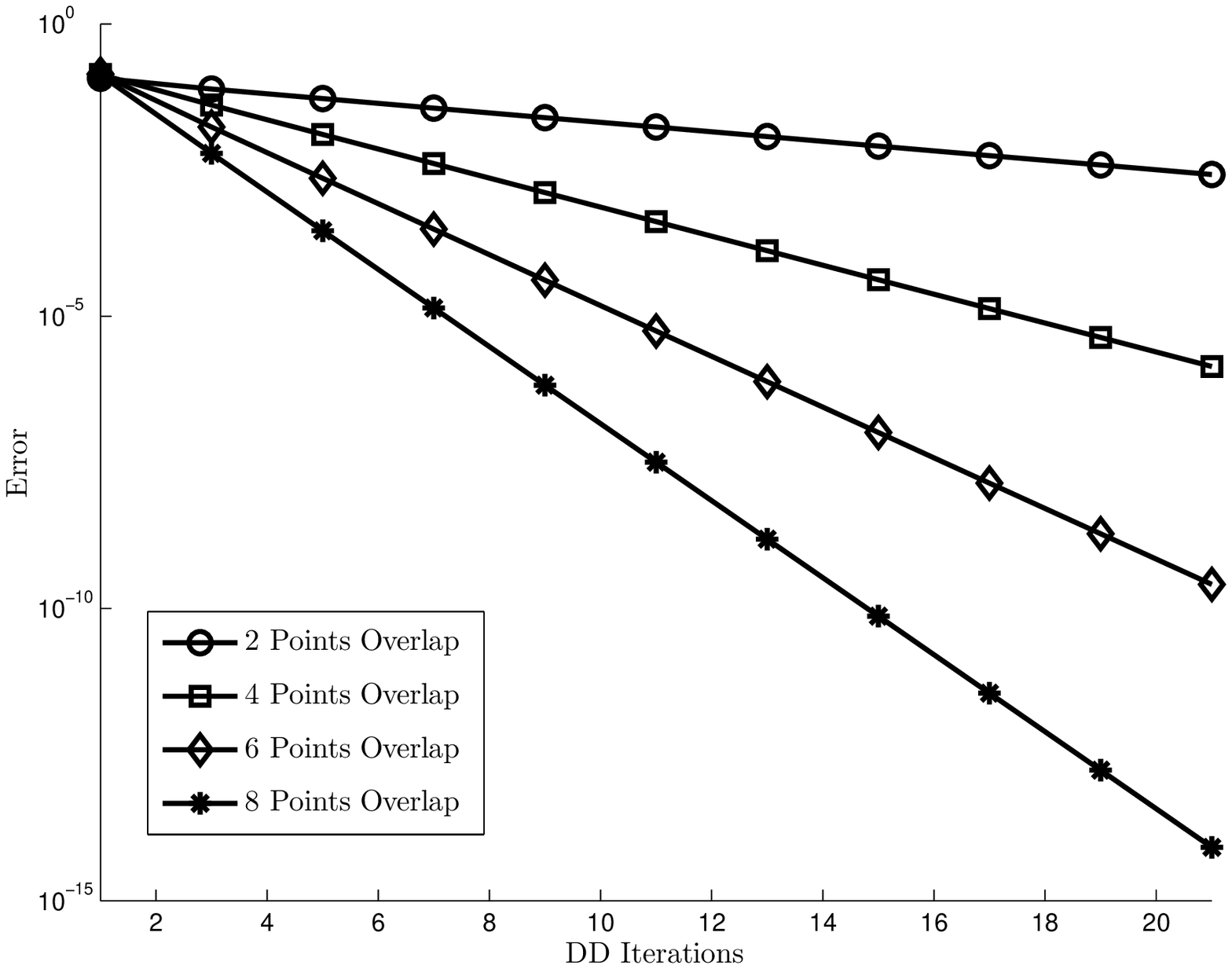}
\caption{Classical Schwarz convergence histories for varying amounts of overlap.}
\label{haynesr_mini_17_figure_DD_2D_comp1_2}
\end{minipage}
\end{figure}

In Figure \ref{haynesr_mini_17_figure_DD_2D_comp1_1} we plot the maximum error between the subdomain and single domain
solutions
for each of $x_1^n$, $x_2^n$, $y_1^n$, and $y_2^n$ obtained 
using classical Schwarz. These are obtained over a 12 by 12 grid with 4 points of overlap in $\xi$ and parameters $a =0.7$ and $b=0.05$. As can be seen, each component of the solution converges at approximately the same rate, so we simplify our discussion by comparing the convergence of only $x_1^n$ in the remaining figures. 

In Figure \ref{haynesr_mini_17_figure_DD_2D_comp1_2} we compare the classical Schwarz algorithm for varying amounts of overlap, 
using $2$, $4$, $6$ and $8$ points of overlap in the $\xi$ direction. As expected, the rate of convergence improves as the overlap increases.

\begin{figure}[h]
\begin{minipage}[b]{0.47\linewidth}
\centering
\includegraphics[width=.95\linewidth]{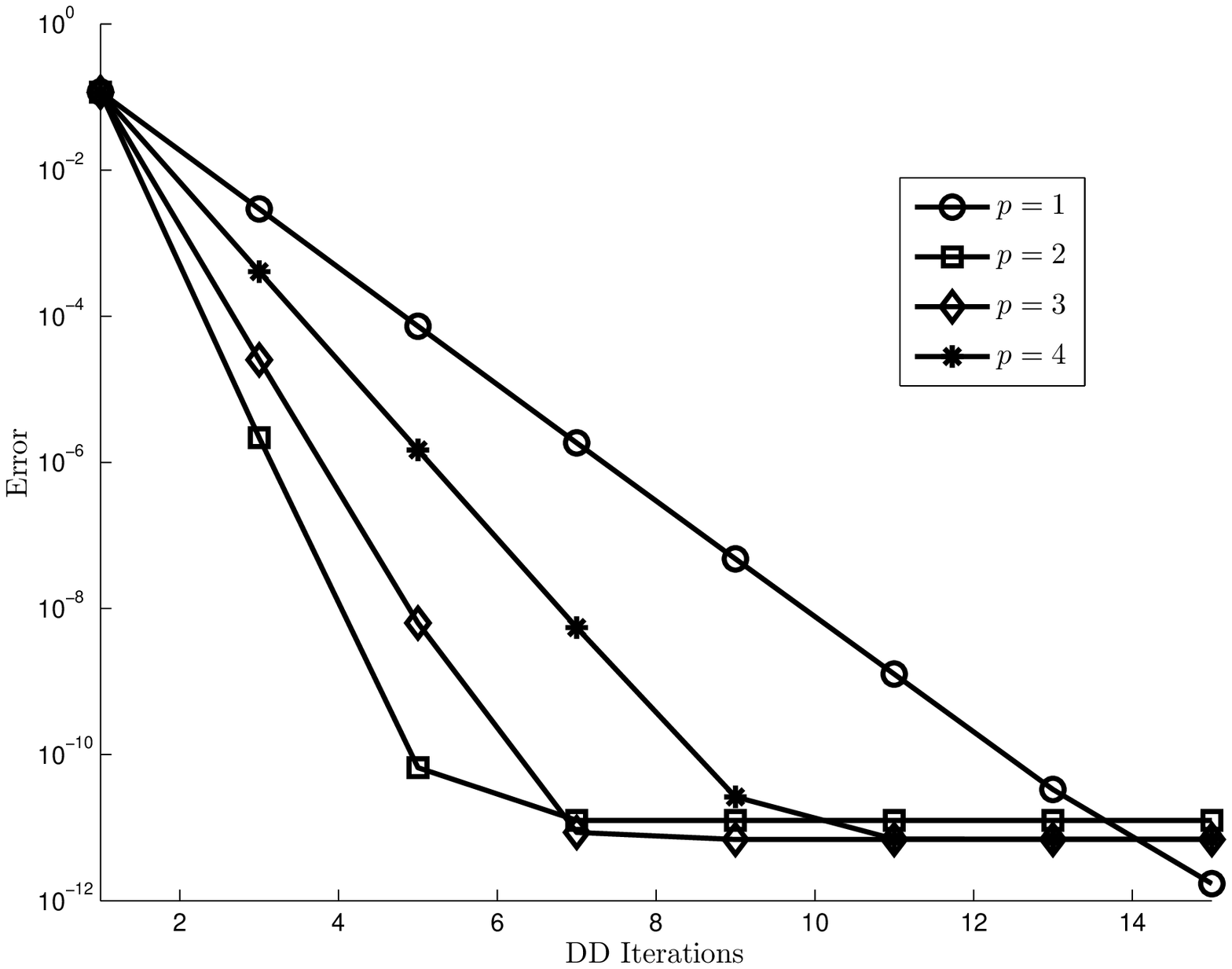}
\caption{Convergence histories for the Schwarz iteration using linear Robin conditions for varying $p$.}
\label{haynesr_mini_17_figure_DD_2D_comp2_1}
\end{minipage}
\hspace{0.5cm}
\begin{minipage}[b]{0.47\linewidth}
\centering
 \includegraphics[width=.95\linewidth]{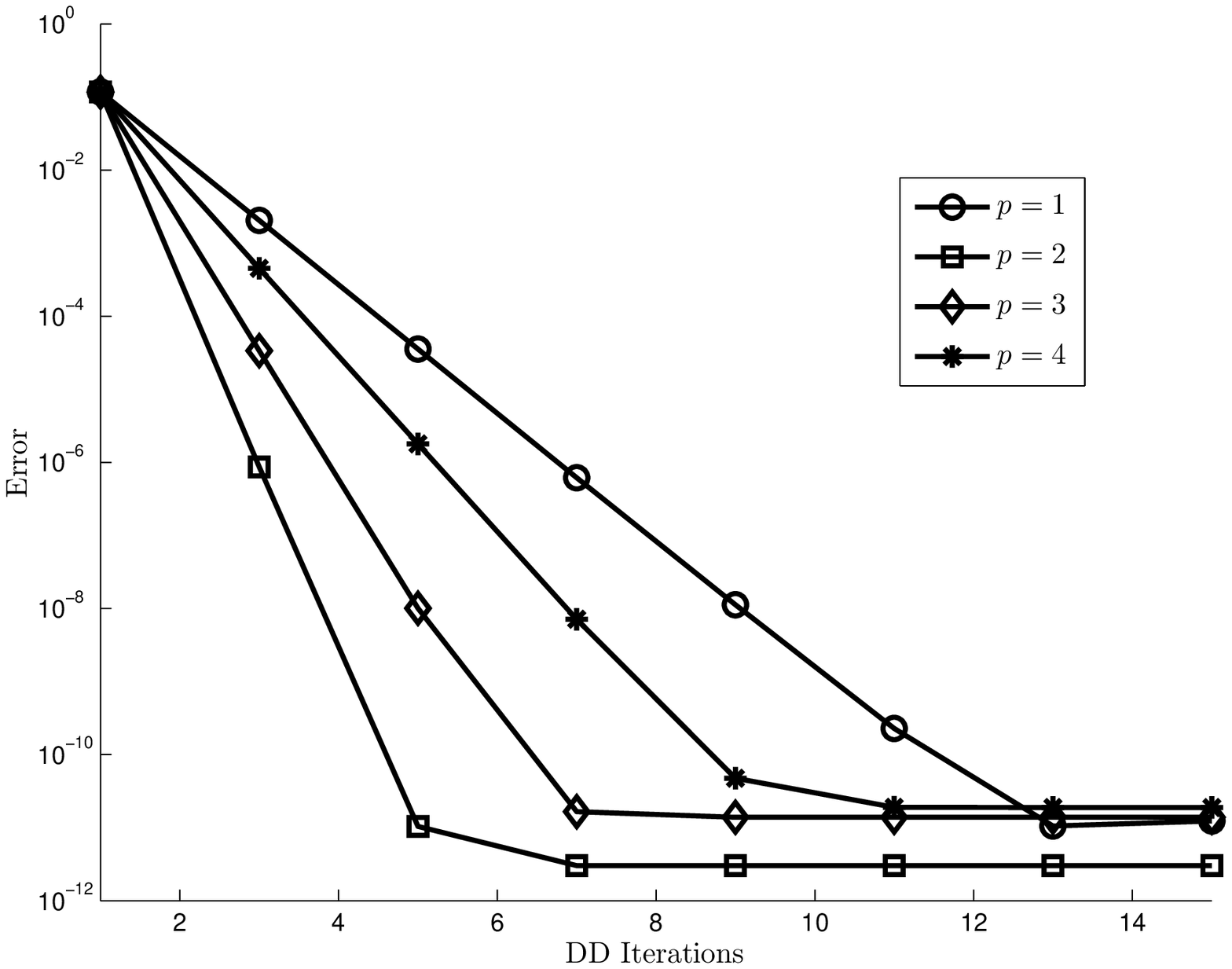}
\caption{Convergence histories for the Schwarz iteration using nonlinear Robin conditions for varying $p$.}
\label{haynesr_mini_17_figure_DD_2D_comp2_2}
\end{minipage}
\end{figure}

\begin{figure}[h]
\begin{minipage}[b]{0.47\linewidth}
\centering
 \includegraphics[width=.95\linewidth]{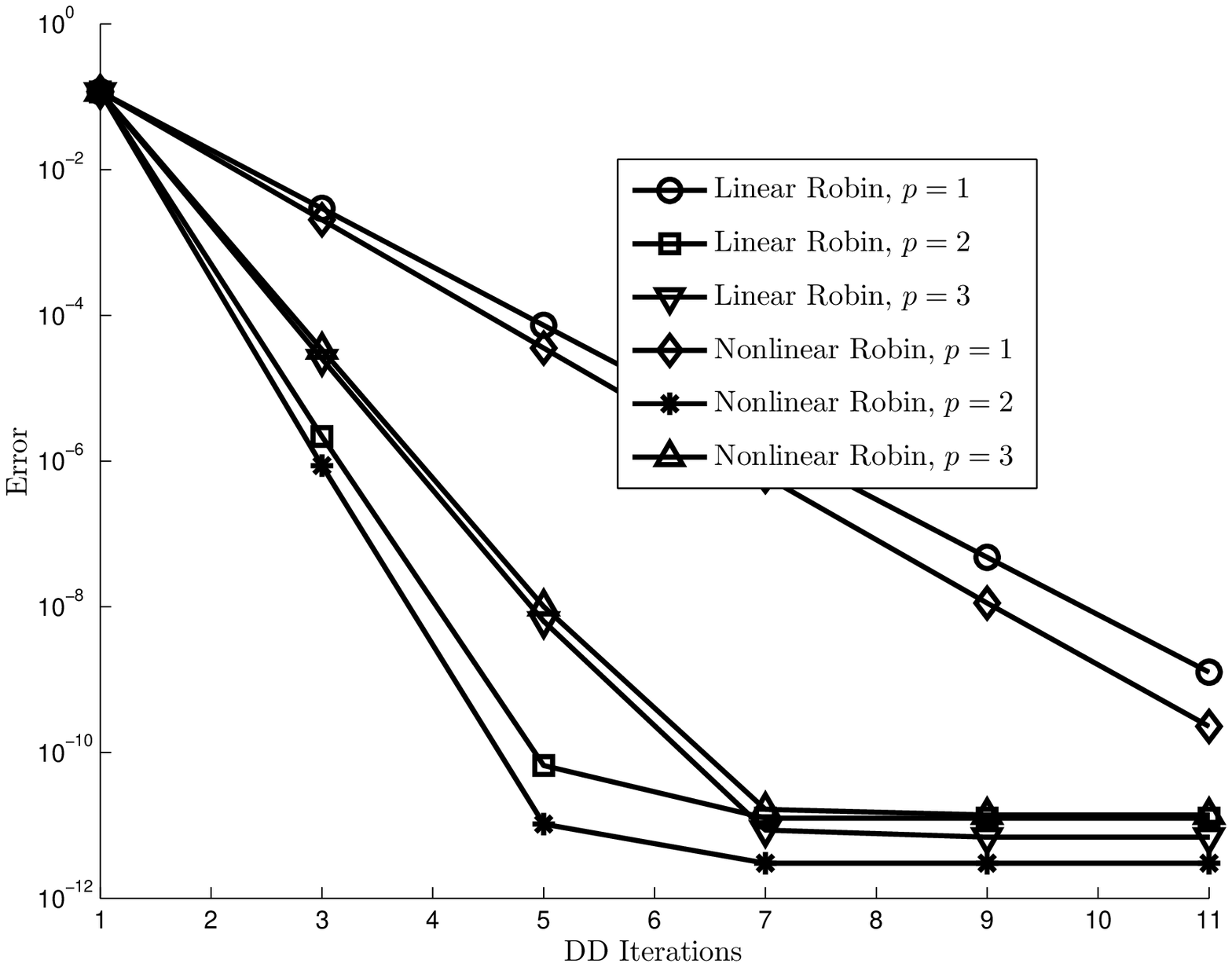}
\caption{Convergence histories for linear and nonlinear Robin transmission conditions with varying $p$.}
\label{haynesr_mini_17_figure_DD_2D_comp3_2}
\end{minipage}
\hspace{0.5cm}
\begin{minipage}[b]{0.47\linewidth}
\centering
\includegraphics[width=.95\linewidth]{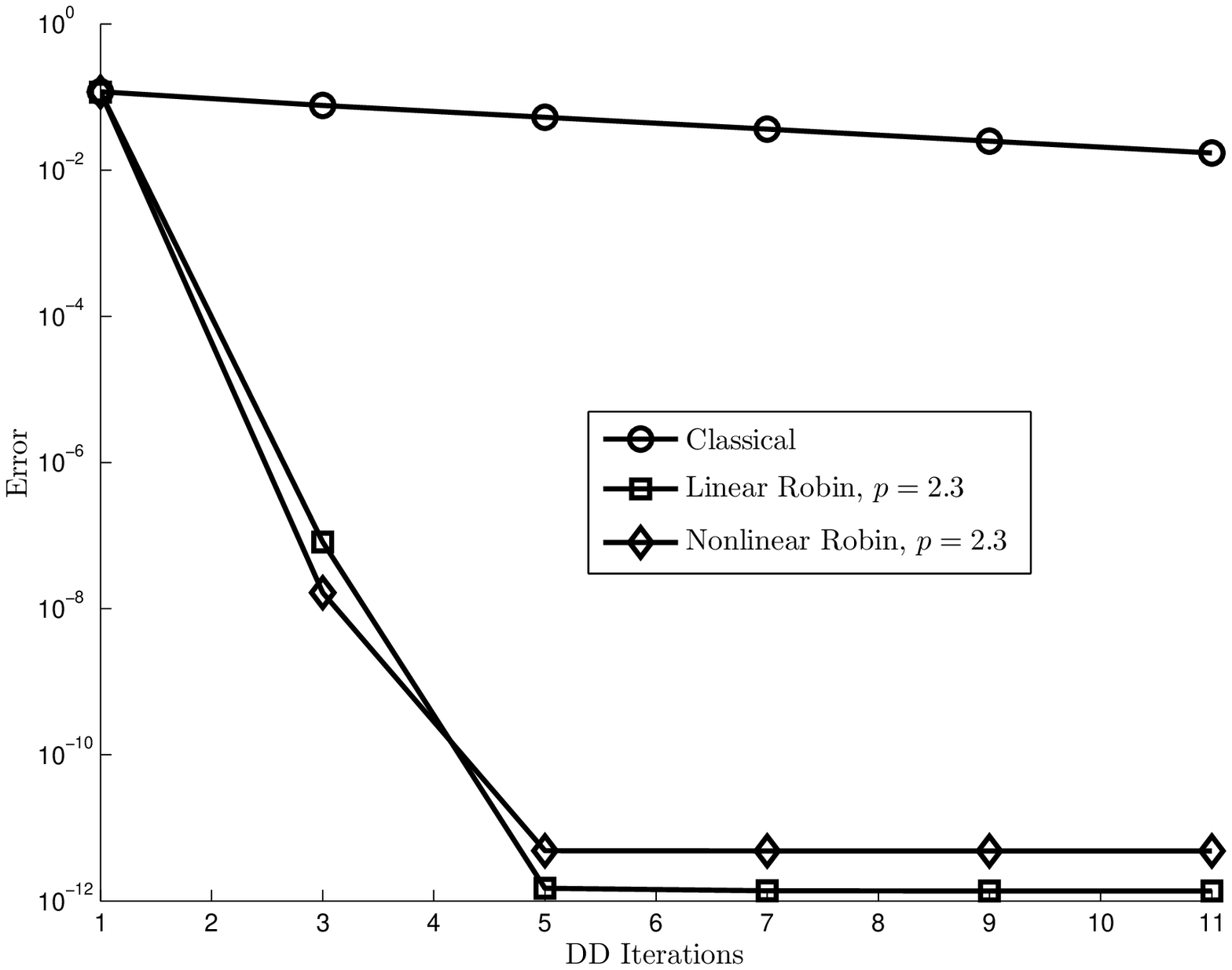}
\caption{Convergence histories for all three iterations considered.}
\label{haynesr_mini_17_figure_DD_2D_comp3_1}
\end{minipage}
\end{figure}

For the two possible optimized Schwarz iterations, we examine the effect of varying the parameter $p$ in Figures \ref{haynesr_mini_17_figure_DD_2D_comp2_1} and \ref{haynesr_mini_17_figure_DD_2D_comp2_2}. To generate these results we use a 12 by 12 mesh with two points of overlap in the $\xi$ direction and parameters $a = 0.7$ and $b = 0.05$. For both types of transmission conditions, the best performance observed occurs for $p=2$. Comparing the linear Robin condition (Figure \ref{haynesr_mini_17_figure_DD_2D_comp2_1}) and nonlinear Robin condition (Figure \ref{haynesr_mini_17_figure_DD_2D_comp2_2}), we see that the convergence histories for a general $p$ are very similar. To examine these similarities, we plot the convergence histories for both optimized iterations for $p=1,2,3$ on the same set of axes in Figure \ref{haynesr_mini_17_figure_DD_2D_comp3_2}. We see that while the variations in this particular case are small, the nonlinear transmission conditions consistently outperform the linear Robin conditions. This is also observed in the results of Figure \ref{haynesr_mini_17_figure_DD_2D_comp3_1}, in which we plot convergence histories for all three proposed DD algorithms. 
For this example we use a 12 by 12 mesh decomposed into two subdomains, with two points of overlap in $\xi$ and parameters $a = 0.7$ and $b = 0.05$. In this example we see that both optimized Schwarz methods vastly outperform classical Schwarz, with the nonlinear transmission conditions slightly outperforming the linear Robin conditions.

Another way to assess the meshes obtained from a DD iteration is to compute a mesh quality measure. 
An equidistribution quality measure for each element $K$ of the grid, $Q_{eq}(K)$, is presented in \cite{haynesr_mini_17_huangbook}.
The  maximum of $Q_{eq}$ over all elements is $1$ if and only if the equidistribution condition is satisfied exactly. 
The larger the value of $max_K Q_{eq}(K)$ the farther the mesh is from equidistributing $\pmb{M}$. 
In Table \ref{haynesr_mini_17_table} we compute the  $max_K Q_{eq}(K)$ for the first five iterations of each proposed Schwarz algorithm. 
The zero column gives the mesh
quality measure for the initial uniform $12 \times 12$ mesh and the $\infty$ column gives the mesh quality measure  for the mesh obtained 
by solving  system (\ref{haynesr_mini_17_2dsystem})
over a single domain.   Note, local equidistribution will not give a value of 1 for the mesh quality measure. 
We see that the meshes obtained by the optimized Schwarz algorithms rapidly give good meshes.
\begin{table}[H]
  \tabcolsep0.3em
  \centering
  \caption{Mesh quality measures for the grids obtained by the proposed Schwarz iterations.}
  \label{haynesr_mini_17_table}
  \begin{tabular}{ c | c c c c c c c }
  Iterations & 0 & 1 & 2 & 3 & 4 & 5 & $\infty$ \\ \hline			
Classical &  1.6375  &  1.3630  &  1.3629  &  1.3178  &  1.3136  &  1.2795  &  1.1979\\ 
Linear Robin & 1.6375  &  2.0076  &  1.1979  &  1.1979  &  1.1979  &  1.1979 &  1.1979\\
Nonlinear Robin & 1.6375  &  2.0114  &  1.1979  &  1.1979  &  1.1979  &  1.1979  &  1.1979\\
  \end{tabular}
\end{table}

\section{Conclusion}
\label{haynesr_mini_17_sec:conclusion}

In summary, we have proposed three different Schwarz DD iterations for obtaining 2D adaptive meshes defined by a local equidistribution principle. The numerical results show that the optimized methods provide a significant improvement over the slow convergence of classical Schwarz, with the nonlinear transmission conditions inspired by the work of \cite{haynesr_mini_17_SINUM} exhibiting the best results. Ongoing work includes the theoretical analysis of these DD approaches for 2D mesh generation and coupling the DD mesh generation with a DD solver for the physical PDE of interest.



\end{document}